# A novel preconditioned conjugate gradient multigrid method for multi-material topology optimization

Nam G. Luu *[1] and Thanh T. Banh [2a]

[1]Faculty of Information Technology, Industrial University of Ho Chi Minh city, Ho Chi Minh city 70000, Vietnam

[2]Department of Architectural Engineering, Sejong University, 209 Neungdong-ro, Gwangjin-gu, Seoul 05006, Republic of Korea

**Abstract.** In recent years, topology optimization has been developed sufficiently and many researchers have concentrated on enhancing to computationally numerical algorithms for computational effectiveness of this method. Along with the development of topology optimization, High Performance Computing (HPC) was marked by a strong dynamic mechanism with a continuous appearance and disappearance of manufacturers, systems, and architectures. Preconditioned conjugate gradient multigrid method (pCGMG) is the most popular in HPC due to its advantage in very large-scale problems. The idea which applies high performance computing to reduce time of running in multi-material topology optimization (MTO) problems with computational time burdens is newly proposed in this article. In multi-material topology optimization procedures, pCGMG is applied for solving linear equation arising from discretization of differential equations. pCGMG is based on mesh size, and then it is powerful to larger scale problems. For the large scale linear static system, minimal compliance-based design is evaluated in this study. This study contributes to a high-performance computing that pCGMG is integrated to an MTO problem, and numerical examples of pCGMG are executed to compare with optimization results in terms of iteration and time-running of different mesh sizes of square wall structure.

**Keywords:** topology optimization; high performance computing; multi-material; preconditioned conjugate gradient, multigrid; numerical methods; parallel computing

## 1. Introduction

Topology optimization is a numerical method that optimizes material layouts within a given design space, for a given set of loads, boundary conditions and constraints with the goal of maximizing the performance of the design. Topology optimization has made an incredible progress as an innovative numerical and design method after the pioneering study by Bendnsøe and Kikuchi (1988). This field has also attracted an immense amount of attention from the researchers. For MATLAB implementation, there are three papers which become always considered by related research. First, Sigmund presented a compact MATLAB implementation of a topology optimization code for compliance minimization of statically loaded structures Sigmund (2001). The total number of MATLAB input lines is 99 including optimizer and finite element subroutine. Second, Jensen and his co-workers published a modified version of Sigmund (2001) that is an efficient topology optimization in MATLAB with 88 lines of code Andreassen (2011). These improvements have been accomplished without sacrificing the readability of the code. Andreassen (2011) may therefore be considered as a valuable successor to the Sigmund (2001) who provides a practical instrument to help for those entering the field of topology optimization. Finally, in 2014, Tavakoli and Mohseni published a code with 115 lines to solve multi-material topology optimization problems, Tavakoli and Mohseni (2014). The alternating active-phase algorithm is based on the splitting of a multi-phase topology optimization problem into a series of binary phase topology optimization sub-problems. In a sequential manner, these sub-problems are solved partially using internal solver which is a traditional binary phase topology optimization solver. In this study, the minimum compliance multi-material topology optimization problem is dealt with and the results are compared to those of Tavakoli and Mohseni (2014).

∗Corresponding author, E-mail: luugiangnam.math@gmail.com

[a] Ph.D. Student, E-mail: btthanhk11@gmail.com

In most topology optimization problems, one of main difficulty is to reduce computational cost and burden, especially in the multi-material topology optimization with large-scale structure cases. Due to the inefficiency of cost and time-running in large-scale problems, there has been many researchers who are interested in reducing the computational cost and time. Bruns *et al.* (2002) considered several algorithms to reduce computationally expensive computing in topology optimization. They modified the basic arc length algorithm and embedded this analysis into the topology optimization problem to enhance the standard arc length method and convergence problems. Colominas and París (2009) presented a parallel computing to solve problems of minimum weight topology optimization of structures with stress constraints formulations. Liu et el. (2011) dealt with a new approach by using guide-weight method to show two formulations of topology optimization with SIMP method. Instead of utilizing SIMP method, Yulin and Xiaoming (2004) and Delgado and Hamdaoui (2019) devoted to the application of the level-set method for topology optimization of basic structure and viscoelastic structures, respectively. Liao *et al.* (2019) presented Triple Acceleration Method (TAM) to accelerate topology optimization problems, which consists of three parts: multilevel mesh, initial value-based Preconditioned Conjugate Gradient (PCG) method, and local-update strategy. TAM accelerates topology optimization in three aspects including the reduction of mesh scale, accelerating solving equations, and decreasing the number of updated finite elements. Although the current work takes SIMP, TAM is extended to other types of topology optimization methods such as level-set-based topology optimization and a Moving Morphable Component (MMC)-based topology optimization (Wang and Benson 2016, Xie *et al.* 2018) or combines with advanced methods such as isogeometric analysis (Wang and Benson 2015a, Wang *et al.* 2015b) and GPU parallel computing (Frutos and Pérez 2017a, Frutos *et al.* 2017b, Xia *et al.* 2017). Wang *et al.* (2020) presented a new high efficiency isogeometric topology optimization to improve mesh scale reduction, solving acceleration and design variables reduction.

In the multi-material problems, Zhang *et al.* (2017) presented a novel approach for multi-material topology optimization based on an MMC-framework. Hoan and Lee (2017, 2019) studied various numerical methods in buckling problems to treat a removal of spurious buckling modes for multiple steel materials and buckling constrained using continuum topology optimization. Banh and Lee (2018a, b) dealt with multi-material problems for thin plate on elastic foundations and Reissner-Mindlin plates by using numerical analysis bases on structures of node. A large-scale multi-material topology optimization for scalable framework 3D with octree-based mesh adaptation was considered by Chin *et al.* (2019). A numerical approach to optimize topologies for structures with the presence of cracks by using multiple materials is presented by Banh and Lee (2018c), while Nguyen *et al.* (2018) continued developing for multiphase carbon fiber reinforcement of concrete structures. A modified leaky ReLU scheme (MLRS) for topology optimization with multiple materials is proposed by Liu *et al.* (2019) and Banh and Lee (2019) presented a numerical method to approach multi-directional variable thickness thin plate with multiple materials topology optimization. The proposed method is general and straightforward, and compatible with existing optimization algorithms.

In the present work, a novel preconditioned conjugate gradient multigrid with V-cycle is proposed to reduce remarkably computational cost to multi-material topology optimization problems and the result is compared to the original code Tavakoli and Mohseni (2014). The goal of this study is to emphasize the effectiveness of multigrid in the field of mean compliance based multi-material topology optimization problems in terms of computational cost and time.

The organization of this study is as follows. In Section 2, multigrid method is presented as well as some numerical method to solve linear equations. In the Section 3, the brief knowledge about topology optimization for multi-material and multi-phase is showed. Section 4 shows computational procedures of the present method. The main focus of this paper is Section 5, where we compare pCGMG V-cycle to several numerical methods in a problem discretized from a partial differential equation (PDE), then results of pCGMG is compared to original code. Finally, conclusions of this study are drawn in Section 6.

## 2. Preconditioned conjugate gradient multigrid (pCGMG) and numerical methods for linear system

### *2.1 Introduction to Multigrid*

Multigrid method in numerical analyses is an algorithm for solving linear equations arising from discretization of differential equations. The main idea of multigrid method is to accelerate the convergence of a basic iterative method by a global correction of the fine grid solution approximation, accomplished by solving a coarse problem. It is known as relaxation, which generally reduces short-wavelength error. In the coarse problem it is similar to the fine grid in that it has also short- and long- wavelength errors. This recursive process is repeated until a grid is reached, where the cost of direct solutions is negligible compared to the cost of one relaxation sweep on the fine grid.

There are three main types of multigrid. They are V-, W- and F-cycle, following five steps which are smoothing – reducing high frequency errors, residual computation – computing residual error after the smoothing operation, restriction – down-sampling the residual error to a coarser grid, prolongation – interpolating a correction computed on a coarser grid into a finer grid, and correction – adding prolongated coarser grid solution onto the finer grid. In more details, we consider linear problem $\mathbf{KU} = \mathbf{F}$, where $\mathbf{K}$ is stiffness matrix from discretization of problem, $\mathbf{F}$ is a vector and $\mathbf{U}$ is a vector solution in region $\Omega$. A typical multigrid V-cycle method uses a sequence of $m+1$ nested discretization grids of increasing fineness: $\Omega_0 \subset \Omega_1 \subset \Omega_2 \subset ... \subset \Omega_m = \Omega$ where $\Omega$ is the finest mesh. Associated with the grid sequence which is a sequence of finite element spaces $V_0 \subset V_1 \subset V_2 \subset ... \subset V_m = V$, we also have a linear system corresponding on each level $n$: $\mathbf{K}_n \mathbf{U}_n = \mathbf{F}_n$ where $\mathbf{K}_n$ is the stiffness matrix obtained from discretization in the mesh $\Omega_n$. In order to transform linear systems on the several grid levels, transfer operators are built between finer and courser grids in the form of linear mappings: restriction operator $R_l : \mathbb{R}^{n_l} \to \mathbb{R}^{n_{l-1}}$ and prolongation operator $P_l : \mathbb{R}^{n_{l-1}} \to \mathbb{R}^{n_l}$ where $n_l$ denotes the number of nodes in grid $\Omega_l$ . The basic type of multigrid is two-level with formula:

$$\mathbf{U}^{\text{new}} = \mathbf{U}_l^{\text{old}} - \mathbf{P}_l \mathbf{K}_{l-1}^{-1} \mathbf{R}_l \left( \mathbf{K}_l \mathbf{U}_l^{\text{old}} - \mathbf{F}_l \right) \tag{1}$$

This multigrid scheme is just one of the possibilities to perform a multigrid method. It belongs to a family of multigrid methods and a so-called multigrid method with $\gamma$ -cycle which has the following compact recursive definition in Algorithm 1.

| |
|---|
| Initialization: Stiffness matrix $\mathbf{K}_l$ , load vector $\mathbf{F}_l$ , restriction operator $\mathbf{R}_l$ and prolongation operator $\mathbf{P}_l$ in the level $l$ .<br>If $l = 0$ then:<br>$\mathbf{U}_l = \mathbf{K}_l^{-1} \mathbf{F}_l$<br>Else:<br>- Pre-smoothing step: Apply the smoother $\vartheta_1$ times to $\mathbf{K}_l \mathbf{U}_l = \mathbf{F}_l$ with the initial guess $\mathbf{U}_l$ .<br>- Restrict to the next coarser grid $\mathbf{F}_{l-1} \leftarrow \mathbf{R}_l \left( \mathbf{F}_l - \mathbf{K}_l \mathbf{U}_l \right)$<br>- Set initial iterate on the next coarser grid $\mathbf{U}_{l-1} = 0$<br>- If $\Omega_l = \Omega$ is the finest grid, set $\gamma = 1$<br>- Call the scheme $\gamma$ times for the next coarser grid<br>$\mathbf{U}_{l-1} \leftarrow \mathrm{MG}^{l-1} \left( \mathbf{U}_{l-1}, \mathbf{F}_{l-1} \right)$<br>- Correct with the prolongated update $\mathbf{U}_l \leftarrow \mathbf{U}_l + \mathbf{P}_l \mathbf{U}_{l-1}$<br>- Post smoothing: Apply the smoother $\vartheta_2$ times to $\mathbf{K}_l \mathbf{U}_l = \mathbf{F}_l$ with the initial guess $\mathbf{U}_l$<br>End |
| Alg. 1 Multigrid $\gamma$ -cycle $\mathrm{MG}_\gamma^l \left( \mathbf{U}_l, \mathbf{F}_l \right)$ |

In practice, only $\gamma = 1$ (V-cycle) and $\gamma = 2$ (W-cycle) are used. The names become clear if one has a look on how they move through the hierarchy of grid.

### *2.2 Preconditioned conjugate gradient*

Conjugate Gradient (CG) method is an algorithm for the numerical solution of particular systems of linear equations whose matrix is symmetric and positive-definite. This method is first presented by Hestenes and Stiefel (1952). The conjugate gradient method is often implemented as an iterative algorithm. It is applicable to sparse systems that are impossible to be handled by a direct implementation or other direct methods. It is a descent method which minimizes the functional $T(u) = \|\mathbf{KU} - \mathbf{F}\|_{\mathbf{K}^{-1}}$ and requires only a single matrix vector multiplication per iteration. In Hestenes and Stiefel (1952) authors proved that the method could solve completely after $n$ times iteration and the convergence rate is given by:

$$\left\|\mathbf{U}-\mathbf{U}_k\right\|_{\mathbf{K}^{-1}} \leq \left\|\mathbf{U}-\mathbf{U}_0\right\|_{\mathbf{K}^{-1}} \left(\frac{\sqrt{\kappa}-1}{\sqrt{\kappa}+1}\right)^k \tag{2}$$

where $\kappa$ is a condition number of the matrix $\mathbf{K}$ and $k$ is an iteration number. There is some notation about the convergence relates condition number $\kappa$. The convergence of linear system depends on the condition number of stiffness matrix. It converges quickly if $\kappa$ is small enough but lager than 1 and very slow if $\kappa$. In case $\kappa \gg 1$, a preconditioner to accelerate convergence of the linear system is used that is a matrix or an operator $\mathbf{M}$ satisfying $\kappa\left(\mathbf{M}^{-1}\mathbf{K}\right) \ll \kappa\left(\mathbf{K}\right)$. The construction of preconditioner $\mathbf{M}$ has to ensure that $\kappa\left(\mathbf{M}^{-1}\mathbf{K}\right)$ should be close to 1 and independent to $n$. Classical preconditioners such as Cholesky factorization, diagonal scaling, or Fractorized Sparse Approximate Inverses (FSAI) cannot provide a mesh independent convergence rate and the preconditioned system is often contrast-dependent as shown in Algorithm 2. Therefore, their utilization may be limited in comparisons with modern, numerically scalable multilevel or multigrid (Amir 2014, Vassilevski 2008).

| |
|---|
| Initialization: Stiffness matrix $\mathbf{A}$, load vector $\mathbf{b}$, preconditioner $\mathbf{M}$, initial value $\mathbf{x}_0$<br>$\mathbf{r}_0 = \mathbf{b} - \mathbf{A}\mathbf{x}_0$<br>For $i = 0,1,\ldots$ do:<br>    $\mathbf{z}_i = \mathbf{M}^{-1}\mathbf{r}_i$<br>    If $\mathbf{r}_i = 0$ then: Stop<br>    Else:<br>        $\beta_{i-1} = \mathbf{z}_i^T\mathbf{r}_i / \mathbf{z}_{i-1}^T\mathbf{r}_{i-1}$, $(\beta_{-1} = 0)$<br>        $\mathbf{p}_i = \mathbf{z}_i + \beta_{i-1}\mathbf{p}_{i-1}$, $(\mathbf{p}_0 = \mathbf{z}_0)$<br>        $\alpha_i = \mathbf{z}_i^T\mathbf{r}_i / \mathbf{p}_i^T\mathbf{A}\mathbf{p}_i$<br>        $\mathbf{x}_{i+1} = \mathbf{x}_i + \alpha_i\mathbf{p}_i$<br>        $\mathbf{r}_{i+1} = \mathbf{r}_i - \alpha_i\mathbf{A}\mathbf{p}_i$<br>    End<br>End |
| Alg. 2 Preconditioned conjugate gradient |

The multigrid $\gamma$-cycle (Algorithm 1) can be utilized as a preconditioner in the preconditioned conjugate gradient. Note that when using a Jacobi smoother, the multigrid algorithm results in a symmetric preconditioner. It is a requirement for being applied in preconditioned conjugate gradient (Vassilevski 2008, Pereira *et al.* 2008). The convergence rate of pCGMG depends on the contrast in the system properties and the development of contrast-independent multilevel preconditioners is an active research topic Amir (2014).

### *2.3 Cholesky method and iteration methods*

The preconditioned conjugate gradient multigrid may become a reasonable method to solve linear systems. To prove this argument by practical works, not only by theoretical works, Cholesky method (or Cholesky decomposition), Gauss-Seidel and Jacobi method are considered to compare with the present pCGMG method. Cholesky method is a decomposition of a Hermitian (or self-adjoint matrix, that is, the element in the *i*-th row and *j*-th column is equal to the complex conjugate of the element in the *i*-th row and *j*-th column), positive-definite matrix into the product of a lower triangular matrix and its conjugate transpose, which is useful for efficient numerical solutions.

This method was first presented by Cholesky for real matrices. Cholesky method decomposes a Hermitian positive-definite matrix $\mathbf{A}$ into $\mathbf{A} = \mathbf{L}\mathbf{L}^*$, where $\mathbf{L}$ is a lower triangular matrix with real and positive diagonal entries, and $\mathbf{L}^*$ denotes the conjugate transpose of $\mathbf{L}$. There is some note for Cholesky decomposition: every Hermitian positive-definite matrix (and thus also every real-valued symmetric positive-definite matrix) has a unique Cholesky decomposition; if the matrix $\mathbf{A}$ is a Hermitian semi-definite positive, the lower triangular matrix $\mathbf{L}$ arises from decomposition allowed to be zeros; $\mathbf{A}$ is symmetric then if $\mathbf{A}$ is positive, $\mathbf{L}$ is also Golub and Van Loan (1983). Algorithm 3 is that of Cholesky decomposition.

```
Initialization: Stiffness matrix A, load vector b
For k = 0 : size(b) do:
    L_kk = sqrt(A_kk)
    For i = k+1 : size(b) do:
        L_ik = A_ik / L_kk
    End
    For j = k+1 : size(b) do:
        For t = j : size(b) do:
            A_tj = A_tj - L_ik L_jk
        End
    End
End
```

Alg. 3 Cholesky decomposition

This algorithm may be briefly interpreted by the following formula:

$$L_{ii} = \sqrt{A_{ii} - \sum_{k=1}^{i-1} L_{ik}^2} \tag{3}$$

$$L_{ij} = \frac{1}{L_{jj}} \sqrt{A_{ij} - \sum_{k=1}^{j-1} L_{ik} L_{jk}}, i > j \tag{4}$$

If $\mathbf{A}$ is symmetric and positive definite, $\mathbf{Ax} = \mathbf{b}$ can be solved by first computing of Cholesky decomposition $\mathbf{A} = \mathbf{LL}^*$, solving $\mathbf{Ly} = \mathbf{b}$ for $\mathbf{y}$ by forward-substitution. Finally $\mathbf{L}^*\mathbf{x} = \mathbf{y}$ for $\mathbf{x}$ is solved by back-substitution. The computational complexity of commonly used algorithms is $\mathrm{O}\left(n^3\right)$, where $n$ is the size of matrix $\mathbf{A}$. For linear systems that can be put into a symmetric form, Cholesky decomposition is an alternative for superior efficiency and numerical stability. The weak point of this method is time-running, especially in large-scale problems and in some very large-scale problems. And out of memory phenomenon may occur during computation processes.

Besides direct methods like Cholesky method, iteration methods such as Jacobi or Gauss-Seidel method are common to solve linear systems. In iteration methods, matrix $\mathbf{A}$ is written as $\mathbf{A} = \mathbf{M} - \mathbf{K}$, where $\mathbf{M}$ is an invertible matrix. Then $\mathbf{x} = \mathbf{M}^{-1}\mathbf{Kx} + \mathbf{M}^{-1}\mathbf{b}$ and iteration formula $\mathbf{x}_{k+1} = \mathbf{M}^{-1}\mathbf{Kx}_k + \mathbf{M}^{-1}\mathbf{b}$ are desired. Jacobi method has this form with $\mathbf{M} = \mathbf{D}$ and $\mathbf{K} = \mathbf{L} + \mathbf{U}$, where $\mathbf{D}$ is diagonal matrix with $\mathbf{D}_{ii} = \mathbf{A}_{ii}$, $\mathbf{U}$ and $\mathbf{L}$ are negative strictly-upper and negative strictly-lower triangular matrix of matrix $\mathbf{A}$. In Gauss-Seidel method, $\mathbf{M} = \mathbf{D} - \mathbf{L}$ and $\mathbf{K} = \mathbf{U}$ are chosen with $\mathbf{D}$, $\mathbf{L}$ and $\mathbf{U}$. In some case, a modified method of Jacobi method named damped Jacobi method which $\mathbf{M} = \omega^{-1}\mathbf{D}$ and $\mathbf{K} = \omega^{-1}\mathbf{D} - \mathbf{A}$ where $\omega \in [0,1]$ is chosen. Iteration formula for damped Jacobi method is $\mathbf{x}_{k+1} = \left(\mathbf{I} - \omega\mathbf{D}^{-1}\mathbf{A}\right)\mathbf{x}_k + \omega\mathbf{D}^{-1}\mathbf{b}$, where $\mathbf{I}$ is identity matrix.

## 3. Multiple materials topology optimization model

Since the pioneering research by Bendsøe and Kikuchi (1988), many approaches have been developed for solving topology optimization problems. Besides topology optimization for single material, multi-material topology optimization is also an attractive problem due to composite structure issues and their structural effectiveness and has made outstanding periods as an innovative numerical and design method. In addition to inherit of main ideas of standard topology optimization, the optimal material distribution of variable densities may be continuously derived. Topology optimization may produce efficient structures with higher stiffness by using additional stiffer materials. It also offers material cost savings in compared with single material structures to arrive at required design performance such as strength and stability. But is produces computational burden due to many materials. In this part, the main idea of multi-phase and multi-material is presented, including numerical interactions to save computational time by using the presented numerical methods in Section 2.

### 3.1 Multi-material interpolation

Minimum compliance based on topology optimization problem is considered for multiple materials within an available design domain $\Omega$ discretized by quadrilateral isoparametric finite element. Following Bendsøe and Sigmund (1999), a domain schematic design of multi-material topology optimization is divided and subdomains are shown in Fig. 1, where $\Omega_{\text{void}}$ ( $j = 1, 2, .., n$ ) is a void material. $\Omega_m^1$ and $\Omega_m^2$ are solid subdomains of material 1 and 2 in a multi-material problem, respectively.

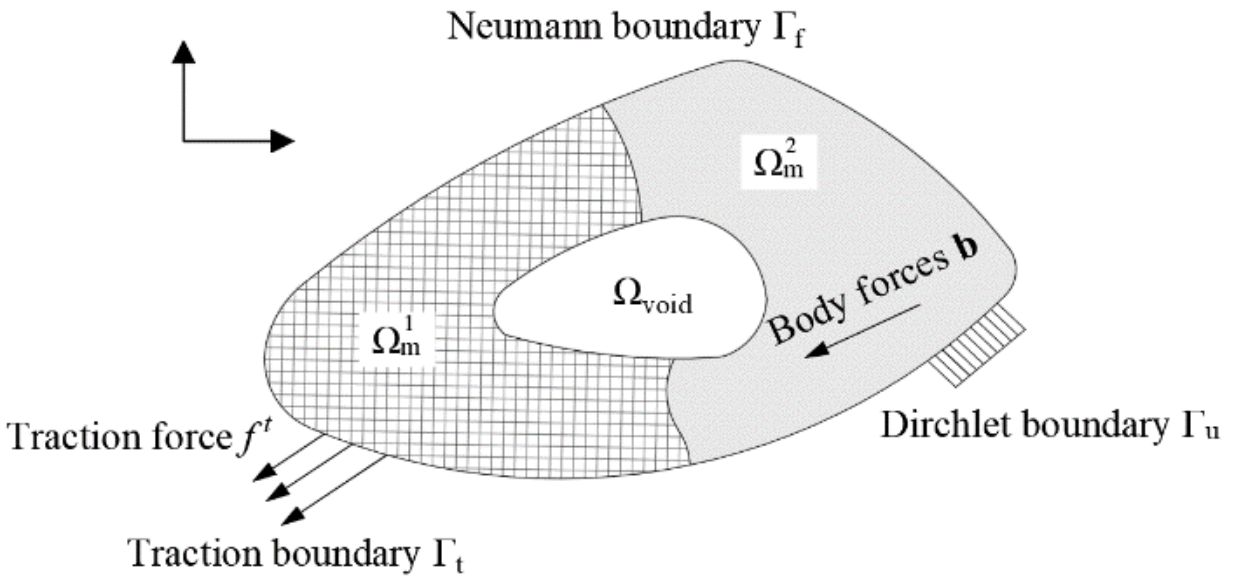

Fig. 1 The design domain, boundary conditions, and external load for multi-material topology optimiz ation

Each design variable is now represented by the density vector $\alpha$ containing element densities $\boldsymbol{\alpha}$. While void is considered as a separate material phase, relative densities of each element are design variables $\alpha$ connected into a vector $\boldsymbol{\alpha}$. Multi-material topology optimization produces the optimal material distribution of $n$ number of materials corresponding to find $n+1$ material phases $\alpha_i = \alpha_i(\mathrm{x})$ at each point $\mathrm{x}$ in the domain $\Omega$. The local stiffness tensor $E$ based on a modified Solid Isotropic Microstructures with Penalization (SIMP) version of a linear interpolation for multiple materials can be formulated by incorporating $\boldsymbol{\alpha}$ as an integer formulation as follows.

$$E(\alpha) = \sum_{i=1}^{n+1} \alpha_i^p E_i^0 \tag{5}$$

where $p$ is the factor that penalizes elements with intermediate densities to approach to almost $0$ or $1$. For this reason, the penalization is reached without presenting any explicit penalization scheme. In Bendsøe and Sigmund (1999), for materials with Poisson's ratio $\nu = 0.3$, authors noted that it is suitable to use $p \geq 3$. $E_i^0$ is elastic modulus corresponding to phase $i$ -th.

### 3.1 Multi-material interpolation

This problem is considered in Tavakoli and Mohseni (2014) as a first example. To avoid singularities in mathematical computing processes of topology optimization, the relaxation for densities between $0$ and $1$ is chosen by a very small lower bound non-zero value $\varepsilon_i$ for this type of problems. The general mathematical formulation of structural multi-material topology optimization problems can be stated as follows:

$$\begin{aligned} &\text{minimize: } C(\alpha_i, \mathbf{U}) = \mathbf{U}^T\mathbf{K}\mathbf{U} \\ &\text{subject to: } \quad \mathbf{K}(\alpha)\mathbf{U} = \mathbf{F} \\ &\qquad\qquad\qquad \int_\Omega \alpha_i d\Omega \le V_i \\ &\qquad\qquad\qquad 0 < \varepsilon_i \le \alpha_i \le 1 \end{aligned} \tag{6}$$

where $C$ is structural compliance, and $\boldsymbol{\alpha}$ is the density vector for phase material $i$-th. $V_i$ is the per-material volume fraction constraint with values of $i$ in range from 1 to $n+1$ with unity $\sum_i V_i = 1$. Global stiffness matrix $\mathbf{K}$, global load vector $\mathbf{F}$ and global displacement vector $\mathbf{U}$ which collects the displacement control variables and additional enrichment degree of freedom can be written as follows.

$$\mathbf{K}_{ij}^e = \int_{\Omega^e} \mathbf{B}_i^u \overline{\mathbf{D}} \mathbf{B}_j^u, \ \mathbf{F}_i = \{\mathbf{F}_i^u\}, \ \mathbf{U} = \{\mathbf{u}\} \tag{7}$$

where $\overline{\mathbf{D}} = \sum_{k=1}^{n+1} \alpha_k^p \mathbf{D}_k^0$ with $\mathbf{D}_k^0$ is the material property matrix corresponding to the phase material including Poisson's ratio $\nu$ and nominal elastic modulus $E_k^0$ is standard strain-displacement matrix as follows.

$$\mathbf{D}_k^0 = \frac{E_k^0}{1-\nu^2} \begin{bmatrix} 1 & \nu & 0 \\ \nu & 1 & 0 \\ 0 & 0 & \dfrac{1-\nu}{2} \end{bmatrix} \tag{8}$$

$$\mathbf{B}_i^u = \begin{bmatrix} N_{i,X_1} & 0 \\ 0 & N_{i,X_2} \\ N_{i,X_2} & N_{i,X_1} \end{bmatrix} \tag{9}$$

For a given materials distribution $\alpha$, the PDE operator $\Re$ could be expressed as follows.

$$
\begin{cases}
\nabla \cdot \left( C : \mathrm{D(u)} \right) & = f(x) \quad \text{in } \Omega \\
u(x) & = u_0(x) \quad \text{on } \Gamma_u \\
\left( C : D(u) \right) \cdot \mathbf{n} & = 0 \quad \text{on } \Gamma_f \\
\left( C : D(u) \right) \cdot \mathbf{n} & = t_0(x) \quad \text{on } \Gamma_t
\end{cases}
\tag{10}
$$

## 4. Computational procedures of multi-material optimization with pCGMG

A briefly computational algorithm of multi-material topology optimization based on pCGMG is shown in Fig. 2. This procedure presents an optimality criteria-based alternating active-phase algorithm using Gauss–Seidel iteration version of multi-material with two active phases `a' and `b', and $`n' > 2$ materials.

In addition, an explicit viewpoint through utilizing pCGMG is considered. To perform the pCGMG step, the geometry, loading and boundary conditions, and material properties are determined. Then, finite element analysis is applied to receive the stiffness matrix $\mathbf{K}$ and load vector $\mathbf{F}$. The linear system $\mathbf{KU} = \mathbf{F}$ is solved by direct solution method and pCGMG in the flowchart. In the right side of the flowchart, the performance of pCGMG is presented in a recursive form. After preconditioning conjugate gradients, the recursion runs through 'level'-1 times where 'level' is fixed before the recursion. Solution $\mathbf{U}$ is the result of the linear system $\mathbf{KU} = \mathbf{F}$ for finite element analysis. By using the solution $\mathbf{U}$, a sensitivity analysis of objective function with respect to element design variables is calculated. The sensitivity filtering is utilized at the same time. The design variables of a binary sub-problem are iteratively updated. The iterative process continues until a desired optimum convergence criterion, such as reaching a given minimum of compliance or reaching a given number of iterations.

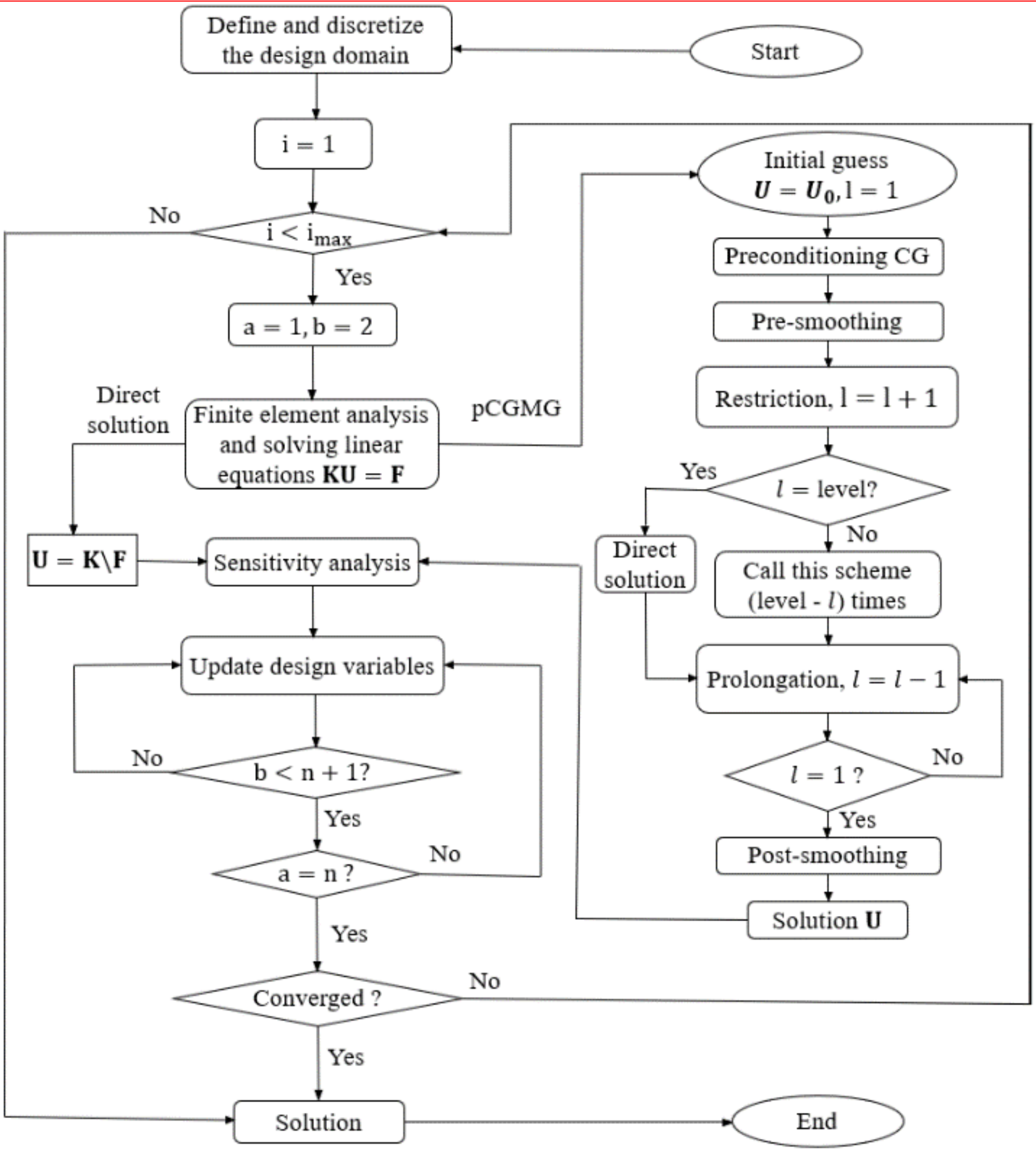

Fig. 2 Flowchart of multi-material topology optimization procedure with preconditioned conjugate gradient multigrid method

## 5. Numerical examples and discussion

### *5.1 Time-running comparisons with pCGMG, Cholesky and iteration methods*

In this section, a partial differential equation named Poisson's equation is considered. This is an elliptic type of PDE to describe the potential field caused by a given charge or mass density distribution. Statement of Poisson's equation is: $\Delta u = f$ in domain $\Omega$ where $\Delta$ is the Laplace operator, and $f$ and $u$ are real or complex-valued functions on $\Omega$. Dirichlet boundary condition, which means $u$ is zeros on $\partial\Omega$. After transferring this equation to weak formulation, finite element and Galerkin method are applied to create stiffness matrix $\mathbf{K}$ and load vector $\mathbf{F}$ for a linear system $\mathbf{KU} = \mathbf{F}$ where $\mathbf{U}$ is a discrete solution of $u$. Multigrid V-cycle, Cholesky, Gauss-Seidel and Jacobi method are dealt with to compare the time-running with five different grids $2^n \times 2^n$ using the same Poisson equation's parameters and maximum level of multigrid. It means direct solution method is utilized in the coarsest mesh $2\times 2$: 1) $16\times 16$ mesh, multigrid level 3; 2) $32\times 32$ mesh, multigrid level 4; 3) $64\times 64$ mesh, multigrid level 5; 4) $128\times 128$ mesh, multigrid level 6; 5) $256\times 256$ mesh, multigrid level 7.

Table 1 The time-running comparison of pCGMG, Cholesky, Gauss-Seidel and Jacobi method in Poisson's equation.

| Method | Mesh | | | | |
|---|---|---|---|---|---|
| | $16\times 16$ | $32\times 32$ | $64\times 64$ | $128\times 128$ | $256\times 256$ |
| pCGMG | 0.23s | 2.62 | 0.56s | 0.65s | 0.41s |
| Cholesky | 0.2s | 4.15s | 151.09s | 7880.85s | *no survey* |
| Gauss-Seidel | 0.027s | 0.03s | 0.25s | 3.66s | 131.05s |
| Jacobi | 0.017s | 0.083s | 0.66s | 10.96s | 523.43s |

According to Table 1, multigrid is the outstanding method to well solve a linear system with large-scale, while Cholesky, Gauss-Seidel and Jacobi only run well in case of small-scales. These three methods produce a long time in large-scale or out of memory (in case Cholesky method). For demonstrating the performance of the standard pCGMG, the results in Amir (2014) are utilized. Convergence of pCGMG depends on the contrast of the stiffness distribution, as expected. It is shown to be mesh-independent which means the same number of pCGMG iteration is required regardless of the resolution of the finest mesh. This is a significant advantage when compared to other preconditioning techniques. Therefore, pCGMG is used to optimally save time-running as well as number of iterations, when multi-material topology optimization is executed.

### *5.2 Minimum compliance topology optimization results using multi-material and pCGMG depending on mesh sizes*

In this section, the performance of the proposed preconditioned conjugate gradient multigrid (pCGMG) method is illustrated. Multi-material topology optimization of a square wall structure as shown in Fig. 3 with three distinct materials is carried out. In order to investigate the time-running effect of pCGMG for multi-material topology optimization, several mesh sizes are applied: $32\times 32$, $32\times 64$, $64\times 64$, $64\times 128$, $128\times 128$ and $128\times 256$. Minimal strain energy is a mean compliance for objective in plain and linear elastostatic states.

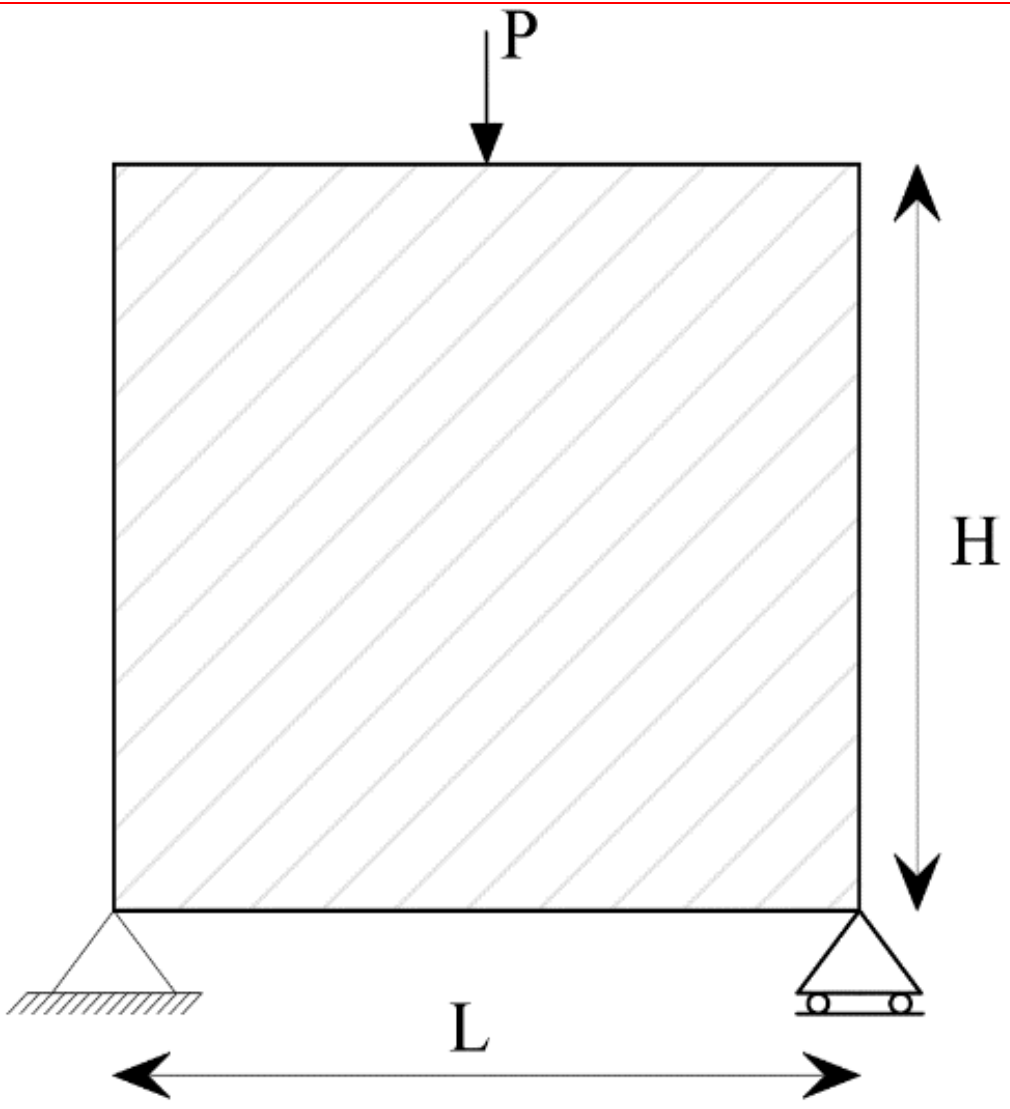

Fig. 3 Geometry, loading and boundary condition of square wall structure depending on mesh size

In Table 2, parameter $p = 4$ illustrates the number of phases considered in the problem while elasticity $\mathbf{e}$, volume fraction are physical properties of all four phases. and denote a filter radius and filter adjusting tolerance, respectively. is a tolerance that is the only solver's stopping criteria, because in this study maximum number of iterations is neglected. In pCGMG's parameters, number of levels is the most important and it is changeable from 2 to 6 to investigate the effect of pCGMG in various levels. The other parameters and are tolerance and maximum number of iterations of preconditioned conjugate gradient, respectively.

Table 2 Design parameters of minimum compliance topology optimization for all numerical methods

| Test problem | $p$ | $\mathbf{e}$ | $\mathbf{v}$ | $rf$ | $tolf$ | $tol$ | $cgtol$ | $cg\,\max$ |
|---|---|---|---|---|---|---|---|---|
| Sq. wall struc. | 4 | $[9;3;1;1e-9]$ | $[0.16;0.08;0.08;0.68]$ | 8 | 0.05 | 1e-3 | 1e-6 | 1000 |

### *5.2.1 Meshes* $32\times 32$ *and* $32\times 64$

By using these meshes, the gap of all time-running results as well as number of iterations are trivial. As can be seen in Table 3, the best choice in terms of time and number of iteration for the mesh $32\times 32$ is pCGMG with $l = 2$ (or called two-level method). With mesh $32\times 64$, the best option based on the number of iteration is pCGMG $l = 3$ and $l = 5$. The difference of the number of iteration to accurately solve is remarkable due to 170 vs 245. However in terms of time-running, the performance of pCGMG of $l = 2$ is better than that of $l = 5$.

Table 3 The comparison of accurate solution and pCGMG in terms of time-running and the number of iterations for minimum compliance topology optimization problem with meshes $32\times32$ and $32\times64$

| Method | Mesh | | | |
|---|---|---|---|---|
| | $32\times32$ | | $32\times64$ | |
| | Iter. | Time | Iter. | Time |
| Accur. solution | 230 | 77.75s | 245 | 174.55s |
| $l=2$ | 207 | 67.91s | 233 | 157.52s |
| $l=3$ | 249 | 77.54s | 170 | 105.36s |
| $l=4$ | 241 | 82.58s | 236 | 147.71s |
| $l=5$ | 248 | 94.45s | 170 | 128.17s |
| $l=6$ | 240 | 95.97s | 188 | 151.54s |

Figs 4a-4b show the change of objective functions of all methods based on iteration steps of meshes $32\times32$ and $32\times64$, respectively. The change of objective function of a pCGMG and direct solution method is negligible in case a mesh $32\times32$. The difference of objective functions in case a mesh $32\times64$ is significantly obvious after 150 times of iteration since pCGMG method reaches the tolerance quickly, while direct solution method needs more nearly 100 times of iteration to stop.

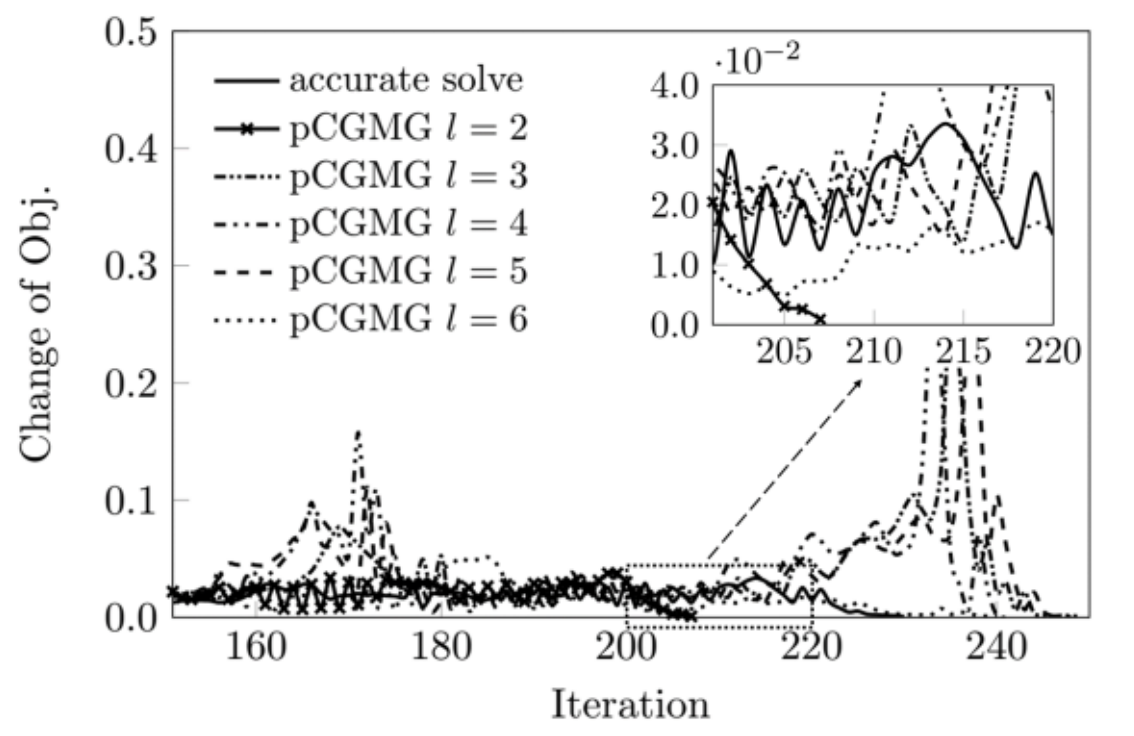

(a) Mesh $32\times32$

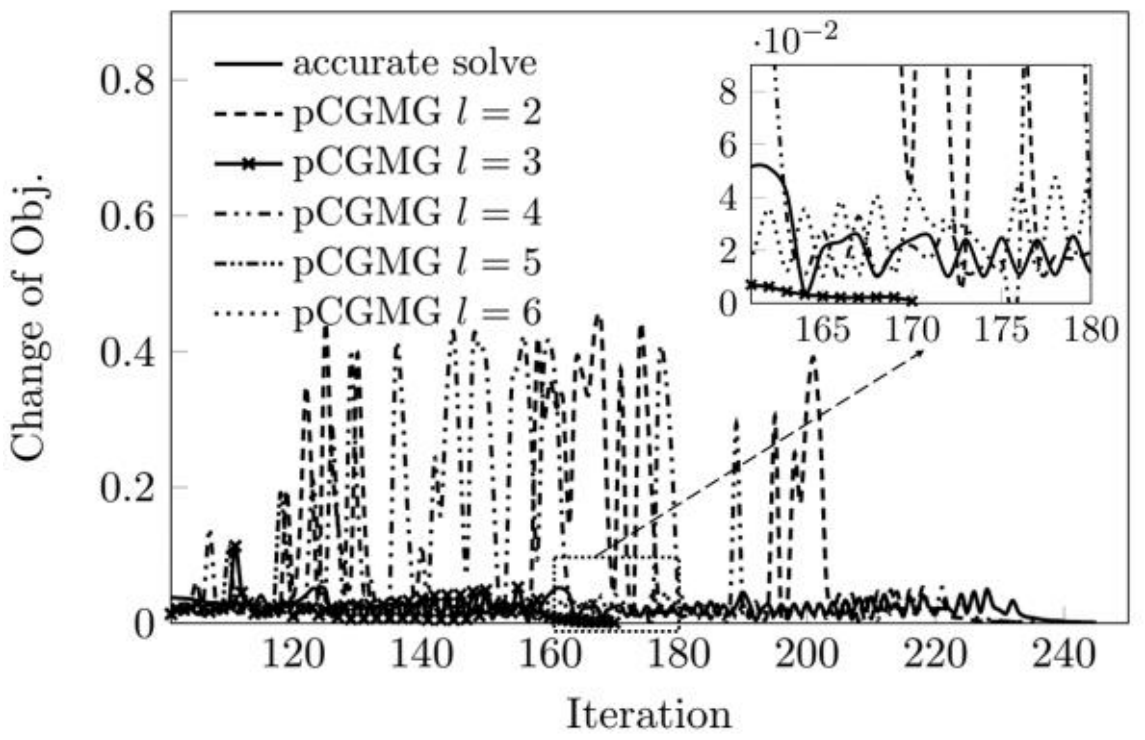

(b) Mesh $32\times64$

Fig. 4 The change of objective function throughout the optimization process of a minimum compliance topology optimization in a square wall structure applied for multigrid methods and direct solution method

Figs 5a-5b illustrate optimal topologies as well as histories of minimum compliance with a use of the best choice of each mesh. Moreover, optimal topologies of two methods are almost the same due to the topology optimization algorithm for the same model. Only linear equation $\mathbf{KU}=\mathbf{F}$ is different each other. Therefore, the proposed method not only reflects the reliability of the optimal solution, but also improves the time-running as well as the variation of error amplitudes as shown in Figs. 4a-4b. Meshes $32\times32$ and $32\times64$ are small meshes, and therefore the converged minimum compliance values of pCGMG with lower level are better than those of higher level.

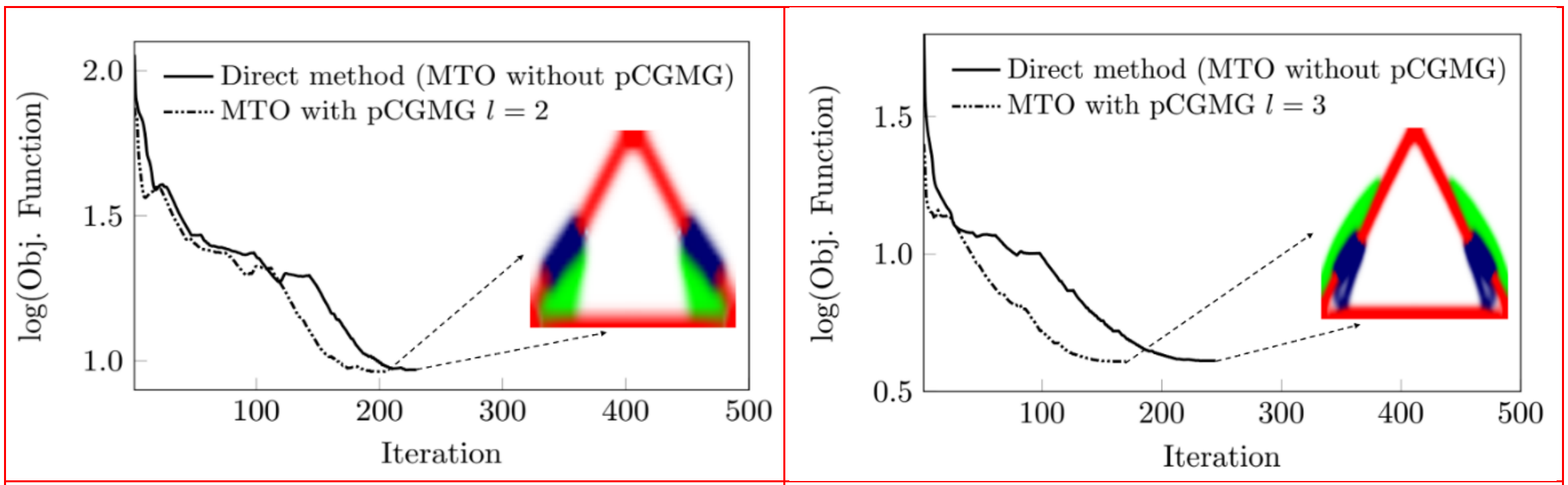

| (a) Mesh $32 \times 32$ | (b) Mesh $32 \times 64$ |
|---|---|

Fig. 5 Histories of objective function and the final optimal topologies with three material problem using pCGMG with the best choice of level.

### 5.2.1 Meshes $32 \times 32$ and $32 \times 64$

As can be seen in Table 4, pCGMG with $l = 3$ is the best choice with respect to time-running and the number of iterations in a mesh $64 \times 64$. In Fig. 6a, oscillations of objective function change in case of pCGMG with level 3 are smaller than that of accurate solutions and converged quickly to the tolerant after 210 times of iteration. In mesh $64 \times 128$, from Table 4, pCGMG method with level 2 and 4 are outstanding in comparing with direct solution methods. However, the time-running of pCGMG with level 4 is shorter than that of level 2. The reason is that pCGMG with level 2 reaches where the cost of direct solution is not negligible in comparing with the cost of one relaxation to sweep on the fine grid. With level 2 of pCGMG method in this case, a mesh $32 \times 64$ is not good for direct solutions after a mesh restriction step.

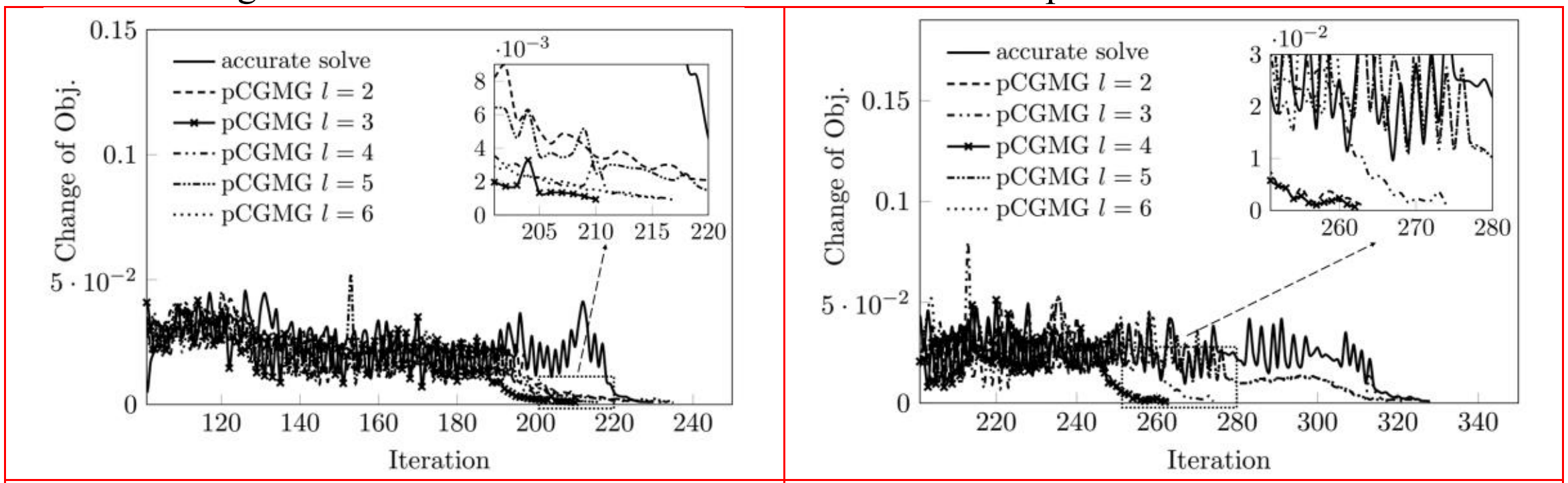

| (a) Mesh $64 \times 64$ | (b) Mesh $64 \times 128$ |
|---|---|

Fig. 6 The change of objective function throughout the optimization process of a minimum compliance topology optimization in a square wall structure applied for multigrid methods and direct solution method

Table 4 The comparison of accurate solution and pCGMG in terms of time-running and the number of iterations for minimum compliance topology optimization problem with meshes $64\times 64$ and $64\times 128$

| Method | Mesh | | | |
|---|---|---|---|---|
| | $64\times 64$ | | $64\times 128$ | |
| | Iter. | Time | Iter. | Time |
| Accur. solution | 230 | 367.08s | 328 | 1079.46s |
| $l=2$ | 232 | 340.07s | 263 | 861.70s |
| $l=3$ | 210 | 281.98s | 274 | 729.43s |
| $l=4$ | 217 | 311.98s | 262 | 1722.01s |
| $l=5$ | 235 | 406.90s | 326 | 987.98s |
| $l=6$ | 216 | 434.78s | 321 | 1124.62s |

Figure 6b demonstrates that objective function changes of pCGMG and direct solution method oscillate constantly until 250-th iteration. After this iteration, pCGMG with converges quickly and stops at 262-nd iteration, while an original code using accurate solution keeps this oscillation and stops nearly 330-th iteration. Figs. 7a-7b describe final optimal topologies as well as histories of compliance of the best choice with meshes $64\times 64$ and $64\times 128$, respectively. With meshes $64\times 64$ and $64\times 128$, the performance of pCGMG is obviously better than that of direct solution method. Level 3 or level 4 of pCGMG is considerable in these meshes.

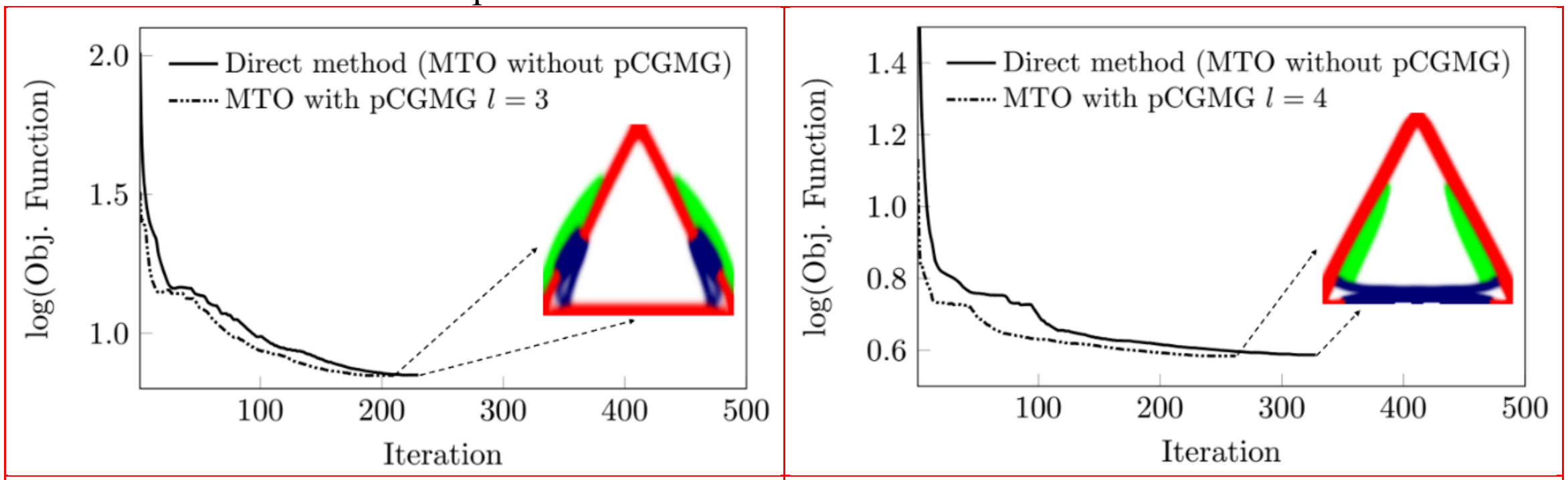

(a) Mesh $64\times 64$ (b) Mesh $64\times 128$

Fig. 7 Histories of objective function and the final optimal topologies with three material problem using pCGMG with the best choice of level.

### *5.2.1 Meshes* $128\times 128$ *and* $128\times 256$

Meshes $128\times 128$ and $128\times 256$ have a large number of finite elements, since the number of nodes are 16,641 and 33,153, respectively. With these meshes, two-level is not good because the direct solution applied after one restriction step is not negligible in compared with the cost of one relaxation to sweep on the fine grid as shown in Table 5. pCGMG with level 3 produces the best result for a mesh $128\times 128$ in both time-running and the number of iteration (270 times and 1524.1s). In a mesh $128\times 256$, pCGMG with level 5 presents the best result and makes a great gap to that of direct solution method (314 times, 4259.65s vs 412 times, 6120.83s).

Table 5 The comparison of accurate solution and pCGMG in terms of time-running and the number of iterations for minimum compliance topology optimization problem with meshes $128\times128$ and $128\times256$

| Method | Mesh | | | |
|---|---|---|---|---|
| | $128\times128$ | | $128\times256$ | |
| | Iter. | Time | Iter. | Time |
| Accur. solution | 401 | 2732.78s | 412 | 6120.83s |
| $l=2$ | 375 | 2780.47s | 319 | 5453.93s |
| $l=3$ | 270 | 1524.10s | 339 | 5097.79s |
| $l=4$ | 322 | 1969.86s | 366 | 4903.15s |
| $l=5$ | 284 | 2388.77s | 314 | 4259.65s |
| $l=6$ | 435 | 3217.66s | 345 | 5197.33s |

The effectiveness of pCGMG is also presented in Figs. 8a-8b. The oscillation of error reduces after each iteration in pCGMG method (level 3 in a mesh $128\times128$ and level 5 in a mesh $128\times256$). The oscillation of direct solution method is constant in a long time and it takes a lot of time to solve. As can be seen, the results demonstrate the stabilities of error of the pCGMG method.

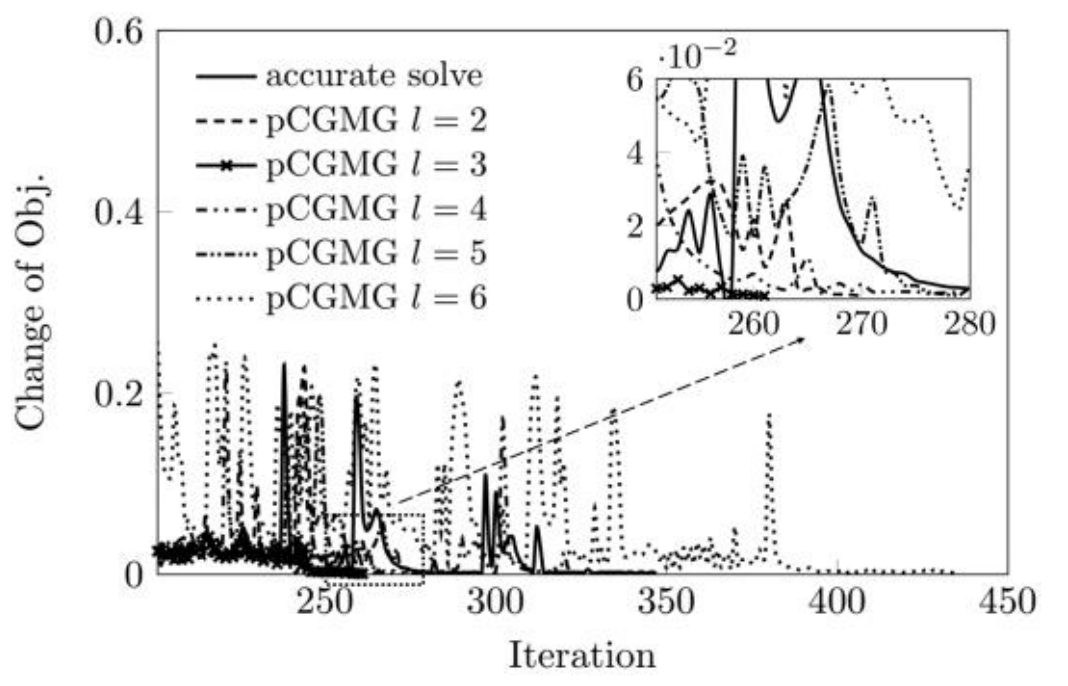

(a) Mesh $128\times128$

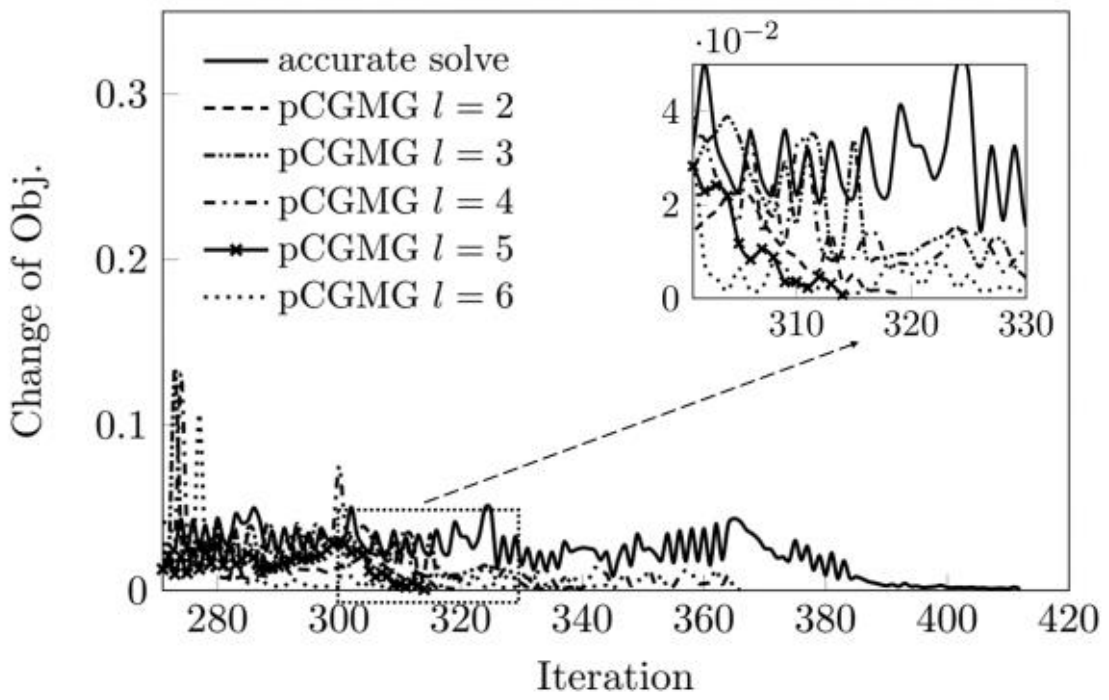

(b) Mesh $128\times256$

Fig. 8 The change of objective function throughout the optimization process of a minimum compliance topology optimization in a square wall structure applied for multigrid methods and direct solution method

Optimal topologies of two methods in these meshes are almost the same as shown in Figs. 9a-9b. Convergence curves of the proposed method reflect the reliability of the optimal solution. Meshes $128\times128$ and $128\times256$ are very large therefore the cost of transferring meshes and relaxation sweeps is negligible to compare with direct solutions of original problems or two-level method. pCGMG methods with level in a range from 3 to 5 are considerable in these cases.

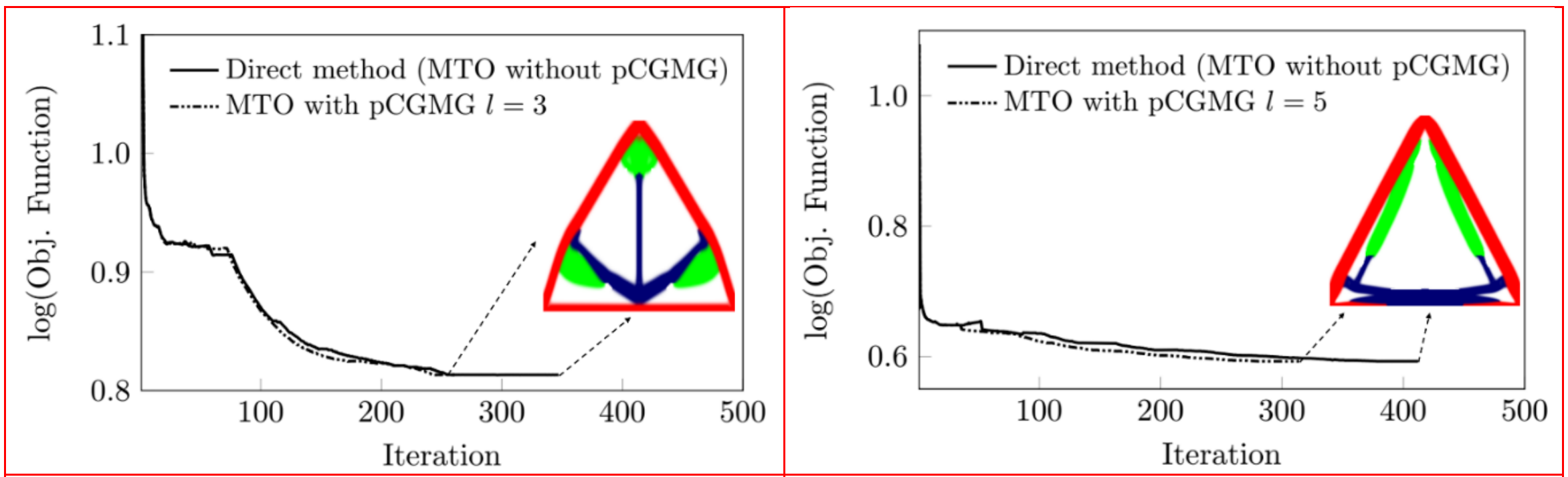

(a) Mesh 128×128 (b) Mesh 128×256

Fig. 9 Histories of objective function and the final optimal topologies with three material problem using pCGMG with the best choice of level.

## 5. Conclusions

This study proposes a novel preconditioned conjugate gradient multigrid method for large-scaled multi-material topology optimization problems. The pCGMG method presents a great efficient computational approach not only on a large linear system raised by discretization of a PDE problem, but also on multi-material topology optimization in large-scale. The effectiveness of pCGMG is to considerably reduce time-running and the number of iteration as well as stabilities of error.

The main finding of this article is to evaluate an appropriate choice of level in pCGMG methods which is oriented on the size of mesh and efficient computation expense. pCGMG with level 2 is the best numerical method in small-scale problems. In large-scale problems, higher levels may be considered because the cost of accurate solution in these levels is negligible in compared with the cost of transferring mesh. The positive results reported here clearly illustrate the suitability of pCGMG methods for applications running in large-scale multi-material topology optimization problems.

For the near future's work, this also encourages a new motivation based on high performance parallel computing environments for enhancing results of crack topology optimization problems or functionally graded materials (FGMs) problems which usually take computational burdens.

## Acknowledgments

This research was supported by a grant (NRF-2017R1A4A1015660) from NRF (National Research Foundation of Korea) funded by MEST (Ministry of Education and Science Technology) of Korean government.